\documentclass[11pt,a4paper]{amsart}
\usepackage[utf8]{inputenc}
\usepackage[english]{babel}
\usepackage{amsmath}
\usepackage{amsfonts}
\usepackage{amssymb}
\usepackage{amsthm}
\usepackage{fullpage}

\usepackage{latexsym}    
\usepackage{amssymb}   
\usepackage{amsmath}    
\usepackage{amsbsy}
\usepackage{amsthm}
\usepackage{amsgen}
\usepackage{amsfonts}
\usepackage{array}
\usepackage[all]{xy}    
\usepackage{tikz}
\usetikzlibrary{cd}
\usetikzlibrary{decorations.pathmorphing}

\usepackage{color}
\usepackage{verbatim} 
\usepackage{url}
\usepackage{enumerate}
\usepackage{enumitem}
\numberwithin{figure}{section}

\tikzset{snake it/.style={decorate, decoration=snake}}

\usepackage{comment}

\DeclareMathOperator{\rank}{rank}

\DeclareMathOperator{\Home}{\mathrm{Hom}}
\DeclareMathOperator{\maxsymrk}{maxsymrk}

\newtheorem{theorem}{Theorem}  
\newtheorem{corollary}[theorem]{Corollary}  
\newtheorem{question}[theorem]{Question} 
\newtheorem{proposition}[theorem]{Proposition} 

		\newtheorem{thm}{Theorem}[section]
		\newtheorem{lem}[thm]{Lemma}

	\theoremstyle{definition}	
		\newtheorem{remark}[thm]{Remark}

\numberwithin{equation}{section}

\newtheorem*{ack}{Acknowledgements}
		
\usepackage{tikz}
\usepackage{subfig}

\usepackage{hyperref}
\hypersetup{
    colorlinks=true,
    linkcolor=blue,
    citecolor=blue,
    urlcolor=blue,
    pdfauthor={Diego Corro and Fernando Galaz-Garcia},
    pdftitle={Positive Ricci curvature on simply-connected manifolds with cohomogeneity-two torus actions}
}

\newcommand{\Z}{\mathbb{Z}} 

\newcommand{\Sp}{\mathbb{S}} 
\newcommand{\CP}{\mathbb{C}P} 

\newcommand{\symrk}{\mathrm{symrk}}

\begin{document}
		
\author[D.~Corro]{Diego Corro$^{*,\dagger}$}
\address[D.~Corro]{Institut f\"ur Algebra und Geometrie, Karlsruher Institut f\"ur Technologie (KIT), Germany.}
\email{\href{mailto:diego.corro@partner.kit.edu}
{diego.corro@partner.kit.edu}}
\urladdr{\url{http://www.math.kit.edu/iag5/~corro/en}}
\thanks{$^*$Supported by CONACYT-DAAD (scholarship number 409912).}


\author[F.~Galaz-Garc\'ia]{Fernando Galaz-Garc\'ia$^{**,\dagger}$}
\address[F.~Galaz-Garc\'ia]{Institut f\"ur Algebra und Geometrie, Karlsruher Institut f\"ur Technologie (KIT), Germany.}
\curraddr{Department of Mathematical Sciences, Durham University, United Kingdom.}
\email{\href{mailto:galazgarcia@kit.edu}{galazgarcia@kit.edu}}
\urladdr{\url{http://www.math.kit.edu/iag5/~galazg/en}}
\thanks{$^{**}$ Supported by the DFG (grant GA 2050 2-1, SPP2026 ``Geometry at Infinity'').}
\thanks{$^{\dagger}$ Supported by the DFG (281869850, RTG 2229 ``Asymptotic Invariants and Limits of Groups and Spaces'').}

\title[Positive Ricci curvature and cohomogeneity-two torus actions]{Positive Ricci curvature on simply-connected manifolds with cohomogeneity-two torus actions}
\date{\today}




	\begin{abstract}
		We show that, for each $n\geqslant 1$, there exist infinitely many spin and non-spin diffeomorphism types of closed, smooth, simply-con\-nect\-ed $(n+4)$-man\-i\-folds with a smooth, effective action of a torus $T^{n+2}$ and a metric of positive Ricci curvature invariant under a $T^{n}$-subgroup of $T^{n+2}$. As an application, we show that every closed, smooth, simply-con\-nect\-ed $5$- and $6$-manifold admitting a smooth, effective torus action of cohomogeneity two supports metrics with positive Ricci curvature invariant under a circle or $T^2$-action, respectively.
	\end{abstract}
	
\maketitle	

\section{Main results}	


	The presence of a smooth (effective) action of a compact (connected) Lie group $G$ on a smooth manifold $M$ has been 
	used to construct Riemannian metrics on $M$ satisfying given geometric properties, e.g.~positive Ricci curvature. In this context, Grove 
	and Ziller showed in \cite{Grove2002} that if $M$ is closed (i.e. compact and without boundary), 
	has finite fundamental group and the action of $G$ has cohomogeneity one (i.e.~the orbit 
	space is one-dimensional), then $M$ admits an invariant Riemannian metric with positive Ricci curvature (see also \cite{ST,Wei1988}). When 	the 
	action is of cohomogeneity two, Searle and Wilhelm showed in 	\cite{Searle2015} that, if the fundamental group of the principal orbits is 
	finite and the orbital distance metric on the orbit space has Ricci curvature greater than or equal to $1$, then $M$ admits an invariant metric 	of 
	positive Ricci curvature. 
	In this article we construct examples of closed, simply-connected smooth manifolds with smooth effective cohomogeneity two 
	torus actions that support Riemannian metrics with positive Ricci curvature invariant under a torus subgroup.

	\begin{theorem}\label{T:MAIN_THM}
		For each integer $n\geqslant 1$, the following hold:
		\begin{enumerate}[label=(\roman*)]
		\item There exist infinitely many diffeomorphism types of closed simply-connected smooth $(n+4)$-manifolds $P$ with a $T^{n}$-invariant Riemannian metric with positive Ricci curvature.
		\item The manifolds $P$ realize infinitely many spin and non-spin diffeomorphism types.
		\item Each manifold $P$ supports a smooth, effective action of a torus $T^{n+2}$ extending the isometric $T^{n}$-action in item (i).
		\end{enumerate}
	\end{theorem}


	Every one of the closed, smooth simply-connected manifolds in Theorem~\ref{T:MAIN_THM}  has a smooth, effective torus action of cohomogeneity 	two. There exist equivariant classification results for these spaces (see \cite{Kim1974},  \cite[Section 4]{Oh1983} and \cite{OR1968}) and,  	in dimensions at most $6$, their topological classification is complete (see \cite{Oh1982,Oh1983,OR1968,Orlik1970}). We combine Theorem~\ref{T:MAIN_THM} with Oh's classification of closed, simply-connected smooth $5$- and $6$-dimensional manifolds with smooth effective cohomogeneity two torus actions (see \cite{Oh1982,Oh1983}) to obtain an explicit list of manifolds admitting metrics with positive Ricci curvature in these dimensions.

	
	\begin{corollary}
	\label{COR:DIM_5_6}
	For every integer $k\geqslant 4$, every connected sum of the form
		\begin{align}
			&\#(k-3)(\Sp^2 \times \Sp^3),\label{trvial5}\\
			&(\Sp^2 \tilde{\times} \Sp^3) \#(k-4) 
						(\Sp^2 \times \Sp^3) \label{nontrivial5},\\
			&\#(k-4)(\Sp^2 \times \Sp^4) \#(k-3) 
						(\Sp^3 \times \Sp^3),\label{trvial6}\\
			&(\Sp^2 \tilde{\times} \Sp^4) \#(k-5) 	
					(\Sp^2 \times \Sp^4)\#(k-3)
						(\Sp^3 \times \Sp^3)\label{nontrivial6},	
		\end{align}			
	has a metric with positive Ricci curvature which is invariant 
	under a torus action of cohomogeneity four. 
	\end{corollary}
	
	Here, $\Sp^2\tilde{\times}\Sp^3$ and 
	$\Sp^2\tilde{\times}\Sp^4$ denote, respectively, the non-trivial $3$-sphere 
	bundle over $\Sp^2$ and the non-trivial $4$-sphere bundle over $\Sp^2$. 
	By convention, we do not consider connected sums where the number of summands of a given type is negative. Thus,  $k\geqslant 5$ in \eqref{nontrivial6}. The manifolds in Corollary~\ref{COR:DIM_5_6} are exactly the closed, simply-connected smooth $5$- and $6$-manifolds that admit effective smooth torus actions of cohomogeneity two (see \cite[Theorem 1.1]{Oh1982} and \cite[Theorem 5.5]{Oh1983}).
	Observe that the manifolds in \eqref{nontrivial5} and \eqref{nontrivial6} are not spin, since their second Stiefel-Whitney class does not vanish (see \cite{Oh1982,Oh1983}). 
	
It is an open question whether every closed, simply-connected smooth $5$-manifold admits a Riemannian metric with positive Ricci curvature. Boyer and Galicki constructed in \cite{Boyer2002,Boyer2006} a countable family of spin, simply-connected  rational homology $5$-spheres which admit a Sasakian metric of positive Ricci curvature. They also exhibited connected sums of such manifolds admitting a Sasakian metric with positive Ricci curvature.  Note, however, that besides the Wu manifold $\mathrm{SU}(3)/\mathrm{SO}(3)$ and the non-trivial bundle $\Sp^2\tilde{\times}\Sp^3$ (which correspond, respectively, to the manifolds $X_{-1}$ and $X_\infty$ in Barden's classification of closed simpy-connected smooth $5$-manifolds \cite{Barden1965}), there are practically no further examples of closed, simply-connected non-spin $5$-manifolds with positive Ricci curvature in the literature. Interestingly, the $5$-manifolds in \eqref{nontrivial5} of Corollary~\ref{COR:DIM_5_6} are non-spin.
	
Without assuming any symmetry,  Sha and Yang showed that connected sums of $k\geq 1$ copies of $\Sp^n\times \Sp^m$, $n,m\geqslant2$, admit metrics of positive Ricci curvature  whose isometry group contains $\mathrm{O}(n)$ (see \cite[Theorem 1 and Remark 1]{Sha1991} and compare with \cite{BGN2003} in the $5$-dimensional case). They also observed that connected sums of finitely many copies of $\Sp^2\times \Sp^4$ and $\Sp^3\times \Sp^3$ also admit metrics with positive Ricci curvature (see \cite[p.\ 129, comments after Theorem 6]{Sha1991}).
Therefore, the manifolds in \eqref{trvial5} and \eqref{trvial6} were already known to admit metrics of positive Ricci curvature. Wraith proved in \cite[Theorem~A]{Wraith2007} the existence of such metrics on connected sums $(\Sp^{n_1}\times\Sp^{m_1})\#(\Sp^{n_2}\times\Sp^{m_2})\#\cdots\#(\Sp^{n_k}\times\Sp^{m_k})$, for $n_i,m_i\geqslant3$ such that $n_i+m_i =n_j + m_j$ for all $1\leqslant i,j \leqslant k$. The spaces in \eqref{nontrivial5} and \eqref{nontrivial6} are, to the best of our knowledge, new examples in the literature, except for the trivial connected sums; the non-trivial bundles over $\Sp^2$  are biquotients (see \cite{DeVito01,DeVito02}) and   Schwachh\"{o}fer and Tuschmann proved that any biquotient with finite fundamental group admits a metric with positive Ricci curvature (see \cite{ST}). We point out that the connected sum decompositions in \eqref{nontrivial5} and \eqref{nontrivial6} are not necessarily equivariant (see \cite{Oh1982,Oh1983}). We also note that, by combining Theorem~\ref{T:MAIN_THM} with \cite[Theorem~0.1]{Gilkey98}, \cite[Main Theorem and results in Section 3]{Nash1979} or \cite[Theorem~0.5]{Wraith1998}, one  can construct further manifolds with positive Ricci curvature.


	 
Closed, simply-connected manifolds with a smooth effective cohomogeneity two torus action are \emph{polar}, i.e.\ there exists an immersed submanifold which intersects all orbits orthogonally (see \cite[Example~4.4]{Grove2012}). Polar manifolds, which include closed cohomogeneity one manifolds, are special in that one can reconstruct them out of the isotropy group information (see \cite{Grove2012}).  It would be interesting to determine whether any other polar manifolds admit (invariant) metrics with positive Ricci curvature.
	 
Although our examples admit smooth effective torus actions of cohomogeneity two, the metrics they support are invariant under a cohomogeneity four action. In the case of cohomogeneity three actions (not necessarily by tori), Wraith has shown that there exist $G$-manifolds with any given number of isolated singular orbits and an invariant metric of positive Ricci curvature; the corresponding result is also true in cohomogeneity five, provided the number of singular orbits is even (see \cite{Wraith2014}). For further background on manifolds with positive Ricci curvature, we refer the reader to \cite{CW}.
	\\

	
We can also consider Theorem~\ref{T:MAIN_THM} in the context of the \emph{symmetry rank}, originally introduced by Grove and Searle in \cite{Grove1994}. Given a closed Riemannian $n$-manifold $M$, its \emph{symmetry rank}, denoted by $\symrk(M)$, is defined as the rank of the isometry group  of $M$. If $M$ has positive sectional curvature, then $\symrk (M) \leqslant \lfloor \frac{n+1}{2} \rfloor$  and the bound is optimal (see \cite{Grove1994}). We call the optimal bound for the symmetry rank of a class $\mathcal{M}^n$ of closed, simply-connected Riemannian $n$-manifolds the \emph{maximal symmetry rank} of the family and denote it by $\maxsymrk (\mathcal{M}^n)$. Considering appropriate Riemannian products of round spheres shows that, in general, the maximal symmetry rank of closed, simply-connected Riemannian $n$-manifolds with positive Ricci curvature is greater than the maximal symmetry rank in the case of positive sectional curvature. The following question then arises naturally:


	\begin{question}
	\label{Q:MSR}
	 What is the maximal symmetry rank for closed, simply-con\-nect\-ed Riemannian $n$-manifolds with positive Ricci curvature?
	\end{question}
	
	It is not difficult to answer this question in dimensions at most six. In dimensions seven and higher, part (i) of Theorem~\ref{T:MAIN_THM}, combined with the fact that no closed, simply-connected smooth $n$-manifold, $n\geqslant 4$, 	admits a cohomogeneity one smooth effective torus action (see, for example, \cite{Galaz-Garcia2010}),  gives further information on the maximal symmetry rank of higher-dimensional Riemannian manifolds with positive Ricci curvature.

	
	\begin{corollary}
	\label{P:MSR}
	Let $\mathcal{M}^n$ be the class of closed, simply-connected Riemannian $n$-manifolds with positive Ricci curvature.
	\begin{enumerate}[label=(\roman*)]
	\item	If $2\leqslant n\leqslant 6$, then $\maxsymrk(\mathcal{M}^n)= \lfloor2n/3\rfloor$.\vspace{2pt}
	\item If $7\leqslant n$, then $n-4\leqslant\maxsymrk (\mathcal{M}^n) \leqslant n-2$.
	\end{enumerate}
	\end{corollary}
	 
Part (i) of Corollary~\ref{P:MSR} also holds if one replaces, in the hypotheses, ``positive Ricci curvature'' with``non-negative sectional curvature'' (see \cite{Galaz-Garcia2010}). In this case, both topological and equivariant classifications for closed, simply-connected Riemannian manifolds with non-negative sectional curvature and maximal symmetry rank are known (see \cite{Galaz-Garcia2014,Galaz-Garcia2010}). Note, however, that in these dimensions the topological and equivariant classifications of closed, simply-connected Riemannian manifolds with positive Ricci curvature and maximal symmetry rank remains open. 
Observe that the upper bound in part (ii) of Corollary~\ref{P:MSR} holds for purely topological reasons. \\

	 
	Recall that a closed smooth manifold $M$ is \emph{string} if it is spin and its first Pontrjagin class satisfies $p_1(M)/2=0$. The \emph{Stolz conjecture} \cite{Stolz1996} asserts that the Witten genus 
	of a closed string manifold $M$ with a Riemannian metric with positive Ricci curvature must vanish (see \cite{Dessai2009} for a survey).
	
	
	\begin{proposition}
	\label{P:STOLZ}
	The Stolz conjecture holds for closed, simply-connected manifolds with an effective torus action of cohomogeneity-two.
	\end{proposition}

	
	Our article is organized as follows. In Section~\ref{S:PROOF_MAIN_THM} we prove Theorem~\ref{T:MAIN_THM}. In Section~\ref{S:PROOF_COROLLARY_B} we prove Corollary~\ref{COR:DIM_5_6}. We prove Corollary~\ref{P:MSR} and Proposition~\ref{P:STOLZ} in Sections~\ref{S:PROOF_PROP_MSR} and \ref{S:PROOF_PROP_STOLZ}, respectively.
		

\begin{ack} The authors would like to thank  Christoph B\"{o}hm, Jonas Hirsch, Wilderich Tuschmann, and David Wraith for helpful conversations on the construction of invariant metrics with positive Ricci curvature, Anand Dessai and Karsten Grove for discussions on the Witten genus, and the referee of a previous version of this article for raising the points which eventually led to the present note. Parts of the work contained in this note were carried out while D.\ Corro was visiting the Department of Mathematics of the University of Notre Dame, and while F.\ Galaz-Garc\'ia was a member of the Mathematical Institute of the University of Bonn and of the Department of Mathematics of Shanghai Jiao Tong University. The authors thank these institutions for their hospitality.
\end{ack}
	

\section{Proof of Theorem~A}
\label{S:PROOF_MAIN_THM}


 We will realize the manifolds in Theorem~\ref{T:MAIN_THM} as total spaces of principal $T^{n}$-bundles over connected sums of finitely many copies of $\pm\CP^2$ or $\Sp^2\times \Sp^2$. We  proceed as follows. We first show that one can choose the base $4$-manifolds and the torus bundles over them so that their total spaces are simply-connected and realize, in each dimension, infinitely many spin and non-spin diffeomorphism types. By the Lifting Theorem of Gilkey, Park and Tuschmann \cite{Gilkey98}, these bundles admit $T^{n}$-invariant metrics of positive  Ricci curvature. Finally, we will show that the total spaces admit smooth effective $T^{n+2}$-actions extending the isometric $T^{n}$ actions.


We start with some topological observations. Although the first two are classical statements, we include their proofs here for the sake of completeness. We denote the  second Stiefel--Whitney class of a smooth manifold $M$ by $w_2(M)$.


\begin{lem}
\label{LEM:SIMPLY_CONNECTED_TOTAL_SPACE}
Let $n\geq 1$.  If $B_k$ is a closed, simply-connected smooth manifold with $H_2(B_k;\Z) \cong \Z^k$, $k\geq n$, then there exists $\alpha\in H^2(B_k;\Z^n)$ such that the smooth principal $T^n$-bundle over $B_k$ with Euler class $\alpha$ has simply-connected total space.
\end{lem}


\begin{proof}
Recall that principal $T^n$-bundles over $B_k$ are classified by the cohomology group $H^2(B_k;\Z^n)$. Since $B_k$ is simply-connected,   the universal coefficient and Hurewicz theorems imply that
\begin{align*}
	H^2(B_k;\Z^n) & \cong \Home_{\Z}(H_2(B_k;\Z),\Z^n) \\
			   &  \cong \Home_{\Z}(\pi_2(B_k),\pi_2(BT^n)).
\end{align*}
By hypothesis, $H_2(B_k,\Z)\cong \Z^k$. Substituting this into the preceding equation, we get  that
\begin{align*}
	 \Home_{\Z}(\Z^k;\Z^n) \cong \Home_{\Z}(\pi_2(B_k),\pi_2(BT^n)).
\end{align*} 
Thus, since $k\geq n$, there exists an onto homomorphism 
\begin{align*}
	\alpha\in \Home_{\Z}(\pi_2(B_k),\pi_2(BT^n)) \cong H^2(B_k;\Z^n).
\end{align*} 
 Let $\pi\colon P\to B_k$ be the  principal $T^n$-bundle over $B_k$ with Euler class $\alpha$. We have the following commutative diagram:
\begin{center}
\begin{tikzcd}
	\pi_2(P)\arrow[r,"{\pi}_{*}"] & \pi_2(B_k)\arrow[r,"\delta"]\arrow{d}[swap]{\alpha} & \pi_1(T^n)\arrow[r,"i_*"] \arrow[d,"\cong"]& \pi_1(P)\arrow[r] & 1\\
	0 = \pi_2(ET^n) \arrow[r] & \pi_2(BT^n)\arrow{r}{\cong}[swap]{\delta}& \pi_1(T^n)\arrow[r,"i_*"] & \pi_1(ET^n)=1 & 
\end{tikzcd}
\end{center}

To see that the central square in the preceding diagram commutes, note that the group homomorphism $\alpha\colon \pi_2(B_k)\to \pi_2(BT^n)$ is induced by a map $f\colon B_k\to BT^n$, i.e.\ $\alpha = f_{*}$ (see, for example, \cite[Lemma~3.1]{Kim1974}). Moreover, by construction,  the bundle $P\to B_k$ is the pullback of the universal bundle $ET^n\to BT^n$ via $f$. Thus the central diagram commutes. Since $\alpha\colon \pi_2(B_k)\to \pi_2(BT^n)$ is onto, the morphism $\delta\colon \pi_2(B_k)\to \pi_1(T^n)$ is onto. Thus $\pi_1(P) = 1$. 
\end{proof}


\begin{lem}\label{L:PULLBACK_SW}
Let $\pi\colon P\to B$ be a smooth principal circle bundle with Euler class $\alpha\in H^2(B;\Z)$ and closed, simply-connected base and total spaces. Then the following hold:
	\begin{enumerate}[label=(\roman*)]
	\item The induced map  $\pi^*\colon H^2(B;\Z_2)\to H^2(P;\Z_2)$ satisfies  $\pi^*(w_2(B))=w_2(P)$.
	\item $w_2(P) =0$ if and only if $w_2(B) = 0$ or $[\alpha] = w_2(B)$ in $H^2(B;\Z_2)$.
	\end{enumerate}
\end{lem}


\begin{proof}
We first prove part (i). The bundle projection map $\pi\colon P\to B$ is a smooth submersion and hence the tangent bundle $TP$ of $P$ splits as
\[
	TP= \pi^\ast(TB) \oplus \varepsilon,
\]
where $\pi^*(TB)$ is the pull-back of the tangent bundle of $B$ and $\varepsilon$ is a real line bundle over $P$. Since $P$ is simply-connected, $\varepsilon$ is trivial. Therefore $w_2(P) = w_2(\pi^\ast(TB)) = \pi^\ast(w_2(B))$.

Now we prove part (ii). The $\mathrm{mod}\ 2$ reduction of the Gysin sequence of the bundle $\pi\colon P\to B$ yields the exact sequence
\begin{align*}
	\cdots \to H^0(B;\Z_2)\stackrel{\smile[\alpha]}{\longrightarrow}H^2(B;\Z_2) \stackrel{\pi^*}{\longrightarrow}H^{2}(P;\Z_2)\to \cdots
\end{align*}
By part (i), $w_2(P)=\pi^*(w_2(B))$. Hence $w_2(P) = 0$ if and only if \linebreak$\pi^*(w_2(B))=0$. By the exact sequence above, this holds if and only if $w_2(B)=0$ or $w_2(B)= \alpha \mod 2$.
\end{proof}


\begin{lem}\label{L:existence_of_non_spin_bundles}
Let $n\geq 1$ and let $B_k$ be a closed, simply-connected smooth $4$-manifold with $H_2(B_k;\Z) \cong \Z^k$. Then the following hold:
\begin{enumerate}[label=(\roman*)]
	\item If $k\geq n$, then there exists a smooth principal $T^n$-bundle $\pi\colon P\to B_k$ with spin simply-connected total space. 
	\item If $k\geq n+1$, then there exists a smooth principal $T^n$-bundle $\eta\colon S\to B_k$ with non-spin simply-connected total space.
\end{enumerate}
\end{lem}

\begin{proof}
We prove the lemma by induction on $n$. The base case $n=1$ is \cite[Theorem~2]{DuanLiang2005}. Suppose now that the statement is true for some $n=m> 1$. Let us prove it for $n=m+1$. We consider first  part (i), the spin case.

Fix $k \geq m+1$ and let $B_k$ be a closed, simply-connected smooth $4$-manifold with $H_2(B_k;\Z) \cong \Z^{k}$. By the induction hypothesis, there exists a principal $T^m$-bundle $\xi\colon Q\to B_k$ such that $Q$ is simply-connected and spin. 
By Lemma~\ref{LEM:SIMPLY_CONNECTED_TOTAL_SPACE}, there exists $\alpha\in H^2(Q;\Z)$ such that the principal circle bundle $\mu\colon P\to Q$ with Euler class $\alpha$ is simply-connected. By Lemma~\ref{L:PULLBACK_SW},
\[
	w_2(P) = \mu^\ast (w_2(Q))
\] 
and $w_2(P)=0$ if and only if $w_2(Q)=0$ or $w_2(Q) = \alpha \mbox{ mod } 2$ in $H^2(Q;\Z_2)$. 
Since $Q$ is spin, $w_2(Q)=0$, so $w_2(P)=0$, i.e.\ $P$ is spin. Since $Q$ is a principal $T^m$-bundle, it has a free action of $T^m$. By \cite[Remarks after Corollary 1.4]{HattoriYoshida1976}, the $T^m$ action lifts to an action on $P$ by bundle isomorphisms. Thus, this lifted action is free on $P$, since the $T^m$ action on $Q$ is free. Furthermore, it commutes with the $T^1$ action on $P$ given by the circle bundle structure (see \cite{Su1963}). Thus we have a principal $T^{m+1}$ bundle $\pi\colon P\to B_k$ with  simply-connected, spin total space. This proves part (i) of the lemma.  

We now consider   part (ii), the non-spin case. Fix $k \geq m+2$ and let $B_k$ be a closed, simply-connected smooth $4$-manifold with $H_2(B_k;\Z) \cong \Z^{k}$. By the induction hypothesis, there exists a principal $T^m$-bundle $\xi\colon Q\to B_k$ such that $Q$ is simply-connected and non-spin. 
The long exact sequence in homotopy for the fibration $\xi\colon Q\to B_k$ yields the following short exact sequence of abelian groups:
\begin{align}
\label{EQ:SHORT_SPLIT}
	0\to \pi_2(Q)\to \pi_2(B_k)\to\pi_1(T^m)\to 1.
\end{align}
Since $\pi_1(T^m)$ is a free abelian group, it is a projective $\Z$-module, so sequence~\eqref{EQ:SHORT_SPLIT} splits. This fact, together with the Hurewicz isomorphism, implies that
\begin{align}
	H_2(B_k;\Z) & \cong \pi_2(B_k) \nonumber\\
			   & \cong \pi_2(Q)\oplus \pi_1(T^m).\label{EQ:DIRECT_SUM_GROUPS}
\end{align}
Since $H_2(B_k,\Z)\cong \Z^k$ and  $k\geq m+2$, it follows from equation~\eqref{EQ:DIRECT_SUM_GROUPS} that $\pi_2(Q)$ is free abelian and $\rank(\pi_2(Q))\geq 2$. Furthermore, by the universal coefficient theorem and the Hurewicz homomorphism, we get that 
\begin{align}
\label{EQ:EULER_CLASS_NON_SPIN}
	 H^2(Q;\Z) & \cong \Home_{\Z}(\pi_2(Q),\pi_2(BT^1)).
\end{align}
Since $\pi_2(BT^1)=\Z$ and $\rank(\pi_2(Q))\geq 2$, there exists an onto homomorphism 
\begin{align*}
\alpha\colon \pi_2(Q)\to\pi_2(BT^1).
\end{align*}
Arguing as in the proof of Lemma~\ref{LEM:SIMPLY_CONNECTED_TOTAL_SPACE}, we see that the principal $T^1$-bundle $\mu\colon P\to Q$ with Euler class $\alpha$ has simply-connected total space.

Let us now see that we can choose the Euler class $\alpha\in H^2(Q;\Z)$ above so that $w_2(Q) \neq [\alpha] $ in $H^2(Q;\Z_2)$, where $[\alpha]$ denotes the map $\alpha \mbox{ mod } 2$.
Recall that $\pi_2(Q)\cong H_2(Q;\Z)$ is free abelian of rank $k-m\geq 2$.
Thus, by the universal coefficient theorem,
\begin{align*}
	H^2(Q;\Z_2) & \cong \Home_\Z(\Z^{k-m},\Z_2).
\end{align*}
Since $k-m\geq 2$, there exists an onto homomorphism $\alpha\colon \Z^{k-m} \to \Z$ such that $[\alpha]\colon \Z^{k-m}\to \Z_2$ is not the second Stiefel--Whitney class $w_2(Q)$.

By Lemma~\ref{L:PULLBACK_SW},
\[
	w_2(P) = \mu^\ast (w_2(Q))
\] 
and $w_2(P)\neq 0$ if and only if $w_2(Q)\neq 0$ and $w_2(Q) \neq [\alpha]$ in $H^2(Q;\Z_2)$. 
Since $Q$ is not spin, $w_2(Q)\neq 0$. By construction, $w_2(Q) \neq [\alpha]$ in $H^2(Q;\Z_2)$. Hence $P$ is not spin. Arguing as in the spin case, we see that $\pi\colon P\to B_k$ is a $T^{m+1}$-principal bundle with non-spin simply-connected total space. This proves part (ii) of the lemma.  
\end{proof}


\begin{remark}
\label{REM:INFINTIE_DIFFEO_TYPES}
In Lemma~\ref{L:existence_of_non_spin_bundles} we have freedom in choosing the rank $k$ of the second homotopy group of the base $B_k$. Then, in each dimension $n\geqslant 5$, we have a countable family of spin and non-spin principal torus bundles with mutually non-diffeomorphic total spaces.
\end{remark}


We are now ready to prove Theorem~\ref{T:MAIN_THM}. 
Let $n\geq 1$, take $k\geq n+1$, and let $B_k$ be a connected sum of finitely many copies of $\pm \CP^2$ or $\Sp^2\times\Sp^2$ so that $H_2(B_k;\Z)\cong \Z^k$. By Lemma~\ref{L:existence_of_non_spin_bundles} and Remark~\ref{REM:INFINTIE_DIFFEO_TYPES}, there exist infinitely many diffeomorphism types of spin and non-spin principal $T^n$-bundles $P_k$ over $B_k$ whose total spaces are simply-connected.

By the work of Sha and Yang \cite{ShaYang1993}, $B_k$ supports a Riemannian metric with positive Ricci curvature. Since $B_k$ is closed and simply-connected, and $P_k$ is a principal $T^{n}$-bundle over $B_k$, the Lifting Theorem of Gilkey, Park and Tuschmann \cite{Gilkey98}  implies that $P_k$ admits a $T^{n}$-invariant Riemannian metric with positive Ricci curvature so that the bundle projection map $\pi\colon P_k\to B_k$    is a Riemannian submersion. Thus, $P_k$ is a closed, simply-connected smooth $(n+4)$-manifold with a $T^n$ invariant metric with positive Ricci curvature.  

We now show that the $P_k$ admit cohomogeneity two torus actions. Orlik and Raymond showed that every closed, smooth simply-connected $4$-manifold with a smooth effective action of $T^2$ is equivariantly diffeomorphic to a connected sum of finitely many copies of $\Sp^4$,  $\pm \CP^2$ or $\Sp^2\times\Sp^2$. It follows then that $B_k$ admits an effective smooth $T^2$-action. By the work of Hattori and Yoshida \cite{HattoriYoshida1976}, any group acting on a closed,  simply-connected smooth manifold $B$ admits a unique commuting lift to any torus bundle $P$ over $B$. Thus, the $T^2$ action on $B_k$ lifts to an action on $P_k$, by bundle isomorphisms. Thus we can define a smooth effective action of $T^{n+2} = T^2\times T^{n}$ on $P_k$ (see \cite{Su1963}). That is, $P_k$ is a closed, smooth simply-connected $(n+4)$-manifold with a smooth, effective action of $T^{n+2}$. By construction, this action extends the isometric action of $T^n$ on $P_k$. This concludes the proof of Theorem~\ref{T:MAIN_THM}. \hfill $\square$


\section{Proof of Corollary~B}
\label{S:PROOF_COROLLARY_B}
The $5$-dimensional case follows from \cite{DuanLiang2005} and the arguments in the proof of Theorem~\ref{T:MAIN_THM}. By \cite[Theorem~2]{DuanLiang2005}, every $5$-manifold listed in \eqref{trvial5} and \eqref{nontrivial5} of Corollary~\ref{COR:DIM_5_6} is a principal circle bundle over a connected sum of finitely many copies of $\pm\CP^2$ or $\Sp^2\times\Sp^2$. Since, by \cite{ShaYang1993}, these connected sums admit metrics with positive Ricci curvature, their total spaces do too, by \cite{Gilkey98}. Note that, as observed in the proof of Theorem~\ref{T:MAIN_THM}, the total spaces of such bundles admit effective smooth $T^3$-actions, so one can also recover their diffeomorphism types using Oh's diffeomorphism classification of closed, simply-connected smooth manifolds with a smooth effective $T^3$-action (see \cite{Oh1983}).

We now consider the $6$-dimensional case. By the constructions in the proof of Theorem~\ref{T:MAIN_THM}, it suffices to show that we can realize each $6$-manifold in \eqref{trvial6} and \eqref{nontrivial6} of Corollary~\ref{COR:DIM_5_6} as the total space of a principal $T^2$-bundle over a connected sum of finitely many copies of $\pm\CP^2$ and $\Sp^2\times\Sp^2$. Observe that the diffeomorphism type of the $6$-manifolds $M$ listed in Corollary~\ref{COR:DIM_5_6} is determined by the integer $k\geq 4$ and whether the second Stiefel--Whitney class $w_2(M)$ is trivial or not. By \cite[Theorem 1.1]{Oh1982}, any closed, simply-connected smooth $6$-manifold with an effective, smooth $T^4$-action must be diffeomorphic to one of the connected sums in \eqref{trvial6} and \eqref{nontrivial6} of Corollary~\ref{COR:DIM_5_6}, with $k\geq4$ corresponding to the number of orbits with isotropy $T^2$. We will use this fact in what follows.

We first show that we can realize the spin diffeomorphism types, i.e.\ the manifolds listed in \eqref{trvial6} of Corollary~\ref{COR:DIM_5_6}. Fix $k\geq 4$, let $m = k-2\geq 2$ and let $B_m$ be a connected sum of finitely many copies of $\pm\CP^2$ or $\Sp^2\times\Sp^2$ with $H^2(B_m;\Z)\cong\Z^m$. Since $m\geq 2$, by part (i) of Lemma~\ref{L:existence_of_non_spin_bundles}, there exists a principal $T^2$-bundle $P$ over $B_m$ with spin simply-connected total space. Moreover, the base $B$ admits a smooth, effective $T^2$-action that lifts to $P$ and commutes with the free $T^2$-action coming from the principal bundle structure of $P$. Thus, we have a closed, simply-connected, smooth spin $6$-manifold $P$ with a smooth, effective torus action of cohomogeneity two. We will now show that the $T^4$-action on $P$ has exactly $k$ orbits with isotropy $T^2$. Let $j$ be the number of such orbits. By the way the $T^4$-action on $P$ was defined, $j$ corresponds to the number of orbits with isotropy $T^2$ of the $T^2$-action on $B_m$. These orbits are isolated fixed points and correspond to the fixed point set $F(B_m)$ of the $T^2$-action on $B_m$ (see, for example,  \cite{Orlik1970}). Let $\chi(\cdot)$ denote the Euler characteristic. It follows from a well-known theorem of Kobayashi (cf.~\cite[Ch.~II, Theorem 5.5]{Kobayashi}) that
\begin{align*}
	 j & =\chi(F(B_m))\\
	   & =\chi(B_m)\\
	   & = 2+m\\
	   & = k.
\end{align*} 
Thus, the $T^4$-action on $P$ has exactly $k$ orbits with isotropy $T^2$. This shows that we can realize all the diffeomorphism types in \eqref{trvial6} of Corollary~\ref{COR:DIM_5_6}. 

We now show that we can realize the non-spin diffeomorphism types, i.e.\ the family \eqref{nontrivial6} in Corollary~\ref{COR:DIM_5_6}. Recall that, for this family, $k\geq 5$. Fix such a $k$, let  $m = k-2\geq 3$ and let $B_m$ be a connected sum of finitely many copies of $\pm\CP^2$ or $\Sp^2\times\Sp^2$ with $H^2(B_m;\Z)\cong\Z^m$. Since $m\geq 3$, by part (ii) of Lemma~\ref{L:existence_of_non_spin_bundles}, there exists a principal $T^2$-bundle $P$ over $B_m$ with non-spin simply-connected total space. The same argument as in the spin case shows that $P$ has a smooth, effective action of $T^4$ with exactly $k$ orbits with isotropy $T^2$.  Thus we can realize all the diffeomorphism types in \eqref{nontrivial6} of Corollary~\ref{COR:DIM_5_6}. This concludes the proof of Corollary~\ref{COR:DIM_5_6}.\hfill $\square$


\begin{remark}
Note that the $6$-manifolds $M$ in the families \eqref{trvial6} and \eqref{nontrivial6} in Corollary~\ref{COR:DIM_5_6} are also determined by the rank of $H_2(M;\Z)$ and whether or not the second Stiefel--Whitney class $w_2(M)$ vanishes. Thus, the constructions in the proof of Theorem~\ref{T:MAIN_THM}, combined with the fact that  the $6$-manifolds in Corollary~\ref{COR:DIM_5_6} are exactly those that admit smooth, effective $T^4$-actions, also show that we can realize every $6$-dimensional diffeomorphism type in Corollary~\ref{COR:DIM_5_6} by means of a principal $T^2$-bundle over some connected sum of finitely many copies of $\pm\CP^2$ and $\Sp^2\times\Sp^2$.

\end{remark}


\section{Proof of Corollary~D}
\label{S:PROOF_PROP_MSR}
Let $M^n$ be a closed, simply-connected Riemannian $n$-manifold with positive Ricci curvature.
Since $M^n$ is simply-connected, $T^n$ cannot act smoothly and effectively on $M^n$. Hence, $\symrk(M^n)\leqslant n-1$. For $n=2$, the Gau\ss--Bonnet theorem implies that $M^2$ must be diffeomorphic to the $2$-sphere. For $n=3$,   $M^3$ must be diffeomorphic to the $3$-sphere, by work of Hamilton \cite{Ha}. The usual round metrics on these spheres have maximal symmetry rank.

For $n\geqslant 4$, no closed, simply-connected smooth $n$-manifold admits a cohomogeneity one smooth effective torus action (see, for example, \cite[proof of Theorem~B]{Galaz-Garcia2010}). Hence, in these dimensions, $\symrk(M^n)\leq n-2$. Considering Riemannian products of round spheres shows that the bound is optimal when $4\leqslant n \leqslant 6$.

Part (ii) follows from part (i) of Theorem~\ref{T:MAIN_THM}, combined with the fact that no closed, simply-connected smooth $n$-manifold, $n\geqslant 4$, admits a cohomogeneity one smooth effective torus action (see, for example, \cite{Galaz-Garcia2010})
\hfill$\square$


\section{Proof of Proposition~E}
\label{S:PROOF_PROP_STOLZ}
	 For $n\geq 3$, an effective smooth action of an $n$-torus $T^n$ on a closed, simply-connected smooth $(n+2)$-manifold  $M$ has no fixed points 	
	 (see, for example, \cite{Oh1983,Raymond1968}). One can therefore find a circle subgroup of $T^n$ that acts without fixed points on $M$. 
	 Hence all Pontrjagin numbers of $M$ vanish (see \cite[Ch.~II, Corollary 6.2]{Kobayashi}). 
	 Since the Witten genus is a linear combination of Pontrjagin numbers, 
	 it follows that the Witten genus of a closed, simply-connected smooth string $(n+2)$-manifold, $n\geqslant 3$, with an effective 
	 action of $T^n$ must vanish. 
	 When $n=2$, a $T^2$-action on a closed, simply-connected smooth $4$-manifold must have fixed points.	
	 However, the Witten genus of a closed, simply-connected smooth string $4$-manifold must also vanish because 
	 the string condition implies the vanishing of the first (and only) Pontrjagin number. Thus, the Stolz conjecture holds for closed, simply-connected manifolds with an effective torus action of cohomogeneity-two.\hfill$\square$



\end{document}